\documentclass[11pt,reqno]{amsart}
\usepackage{amssymb,amsfonts,amsxtra,amsmath,enumerate,tikz-cd,xcolor,mathrsfs,verbatim,centernot,comment}
\usepackage[all]{xy}

\newtheorem{theorem}{Theorem}[section]
\newtheorem{lemma}[theorem]{Lemma}
\newtheorem{corol}[theorem]{Corollary}
\newtheorem{proposition}[theorem]{Proposition}

\newtheorem{definition}[theorem]{Definition}
\newtheorem{remark}[theorem]{Remark}
\newtheorem{example}[theorem]{Example}

\newtheorem*{mainthm}{Main Theorem}

\title{Chern classes of vanishing cycles for maps having no blowups}
\author{Xia Liao}
\address{
Department of Mathematical Sciences, 
Huaqiao University, 
Chenghua North Road 269, 
Quanzhou, Fujian, China
}
\email{xliao@hqu.edu.cn}

\begin{document}

\begin{abstract}
We study the nearby and vanishing cycles for maps with no blowups. For holomorphic maps between complex manifolds admitting at most ICIS, we show that the Riemann-Hurwitz defect in Chern classes is compensated for by the Chern class transformation of the Milnor number constructible function.
\end{abstract}

\maketitle

\section{introduction}

It is safe to say that the Riemann–Hurwitz formula is among the first theorems students of algebraic geometry learn, due to its simplicity, beauty and usefulness. One major result of this paper is a generalisation of the Riemann-Hurwitz formula for maps between complex manifolds admitting at most isolated complete intersection singularities (ICIS). 

Recall that a holomorphic map germ $g:(\mathbb{C}^{m},0) \to (\mathbb{C}^n,0)$ is said to define an ICIS of dimension $m-n$ at $0\in \mathbb{C}^m$ if the $n$ coordinate functions $g_1,\cdots,g_n$ form a regular sequence in $\mathscr{O}_{\mathbb{C}^m,0}$ and the punctured fibre $(g^{-1}(0)\setminus 0,0)$ is smooth. This is equivalent to any of the following two conditions:
\begin{enumerate}
\item $\mathscr{O}_{\mathbb{C}^m,0}$ is a flat $\mathscr{O}_{\mathbb{C}^n,0}$-module via $g^*$ and $(g^{-1}(0)\setminus 0,0)$ is smooth.
\item $m\geq n$ and $g$ is of finite singularity type \cite[Exercise 2 on p.630]{MR568901}.
\end{enumerate}

According to Hamm \cite{MR286143} (see also \cite{MR3112993}), there is a local Milnor fibration at $0\in \mathbb{C}^m$ with smooth fibre $F_0$ homotopic to a bouquet of spheres. Namely,
\[
F_0 \simeq S^{m-n}\vee \cdots \vee S^{m-n}.
\]
If $\mu_g$ is the number of the $(m-n)$-dimensional spheres appearing in the bouquet and $\chi_g$ is the topological Euler characteristic of $F_0$, then 
\[
\chi_g=1+(-1)^{m-n}\mu_g.
\]

Consider now a holomorphic map $f:M\to N$ between complex manifolds, with $\dim M=m$ and $\dim N=n$. Assume that at any $x\in M$ the germ $f_x:(M,x) \to (N,f(x))$ defines an ICIS of dimension $l=m-n$ so that the local Milnor fibration at $x$ gives rises to the numbers $\chi_f(x)$ and $\mu_f(x)$. As functions on $M$, $\chi_f$ and $\mu_f$ are constructible (see \S \ref{section:CSM} and Corollary \ref{coro:constructibility}). Moreover, let
\[
\bar{c}_i(f)=c_i(TM-f^*TN)
\]
be the virtual Chern classes of the map $f$. Let $c_*(\mu_f)$ be the Chern class transformation, or the Chern-Schwartz-MacPherson class (CSM class) of the Milnor number constructible functions $\mu_f: M \to \mathbb{Z}$. We will briefly review the neccessary background of the Chern (CSM) class in \S\ref{section:CSM}. Our main theorem is a unversal expression of the CSM class of $\mu_f$:

\begin{mainthm}
Let $f:M\to N$ be a flat holomorphic map between complex manifolds of relative dimension $l=m-n\geq 0$ whose fibres admit at most isolated singularities. Then
\begin{align}
(-1)^lc_*(\mu_f)&=\Big(c(f^*TN)\cdot (1+\bar{c}_1(f)+\cdots + \bar{c}_{l}(f))-c(TM)\Big)\cap [M] \tag{$\ast$}\\
&=\Big(-c(f^*TN)\cdot\sum_{i=l+1}^{\infty}\bar{c}_{i}(f)\Big)\cap [M]. \notag
\end{align}
\end{mainthm}

\begin{remark}
In fact, this formula was conjectured by Toru Ohmoto two decades ago in his study on higher Thom polynmomials for singularities of maps based on his theory of equivariant CSM classes, see \cite{MR2197579,MR3585782,MR5040859} for details. In particular, the higher Thom polynomials of the function $\mu_f$ is the universal expression of the Segre-Schwartz-MacPherson class as a power series of virtual Chern classes asociated to the map:
\[
tp^{\rm SM}(\mu)=(-1)^{l+1}(1+c_1+c_2+\cdots)(\bar{c}_{l+1}+\bar{c}_{l+2}+\cdots).
\]
The above theorem confirms this expression.
\end{remark}

Since
\[
\bar{c}_{i}(f)=0
\]
for $i\geq l+1$ when $f$ is submersive everywhere, the intuive meaning of the formula's right hand side is the failure of the map $f$ to be everywhere submersive (or a $C^\infty$ fiber bundle in case $f$ is proper) measured in Chern classes of the source and the target. Loosely speaking, the term
\[
c(f^*TN)\cdot (1+\bar{c}_1(f)+\cdots + \bar{c}_{l}(f))
\]
can be thought of as the product of the Chern classes of the base and the generic fibre. The equality then expresses the failure by the Chern class transformation of the Milnor number constructible function.

\begin{corol}[\cite{MR1819638}]
If both $M$ and $N$ are furthermore compact, then
\[
(-1)^l\int \mu_f=\chi(N)\chi(F)-\chi(M)
\]
where $F$ is a generic fibre and $\chi$ is the topological Euler characteristic.
\end{corol}

\begin{proof}
Taking the degree of the equation $(\ast)$, we get
\[
(-1)^l\int \mu_f=c_n(TN)\cap f_*(\bar{c}_l(f)\cap[M])-\chi(M).
\]
Let $f_*(\bar{c}_l(f)\cap[M])=d[N]$. Restricting this equality to an open subset of $N$ over which $f$ is a trivial bundle with fibre $F$, we get $d=\chi(F)$. The corollary follows immediately.
\end{proof}

When $m=n=1$ this corollary becomes the classical Riemann-Hurwitz formula, because $\chi(F)=\deg f$ and at any $x\in M$ we have $\mu_f(x)=e_x-1$ where $e_x$ is the ramification index. For this reason, I also call the right hand side of equation $(\ast)$ the Riemann-Hurwitz defect in Chern classes.

The most essential step of the proof of the main theorem is the construction of the characteristic cycle of $\chi_f$ in $T^*(M\times N)$, where $\chi_f$ is regarded as a constructible function on $M\times N$ via the graph embedding $M\to M\times N$. The construction is best understood in the larger context where we allow singular spaces to be the source of $f$. From now on let $Z$ be a closed irreducible analytic subset of $M\times N$ and let $f: Z\to N$ be the map induced by the projection $M\times N \to N$. Assume $f$ is open and has no blowup in an appropriate sense (see Definition \ref{definition:noblowup}). In the new setup we do not assume $m\geq n$.

\begin{equation*}
\begin{tikzcd}
Z \arrow[rd,"f"] \arrow[r] & M \times N \arrow[d] \\
                                  & N
\end{tikzcd}
\end{equation*}

We need the no blowup condition to define a constructible function $\textup{Eu}^\vee_f$ on $Z$. This function was originally introduced by Sabbah in \cite{MR804052} when $\dim N=1$, and subsequently in \cite{S1984} under other conditions. When $Z$ is the graph of a flat holomorphic map $f:M\to N$ admitting at most ICIS, we have $\chi_f=(-1)^{m-n}\textup{Eu}^\vee_f$. To explain what $\textup{Eu}^\vee_f$ is, let $Z^{\circ}$ be the Zariski open subset of $Z$ where $Z$ is smooth and the rank of $df$ is $n$. We associate to $f$ the relative conormal space $T^*_fM$, which is the closure in $T^*M \times N$ of the set
\begin{equation*}
\{(x,\alpha,y)\ |\ (x,y)\in Z^{\circ} \text{ and } \alpha\in T^*_xM \text{ is conormal to the fibre } f^{-1}(f(x,y)) \subset M\}.
\end{equation*} 
Define $\textup{Eu}^{\vee}_{f}(z)$ for $z\in Z$ as the local intersection number in $T^*(M\times N)$ of the graph of the differential of a local distance function centred at $z$ with the cycle $[T^*_fM]$, corrected by the sign $(-1)^n$. This is essentially the same definition given in \cite{S1984}, but we choose to formulate it slightly differently in order to emphasise the analogy to a well-known construction in the theory of charcteristic cycles which produces a constructible function from a conical Lagrangian cycle by means of a similar intersection. Unless in very special cases $[T^*_fM]$ is not conical Lagrangian, but the crucial point is that we can deform $[T^*_fM]$ to a conical Lagrangian cycle by finitely many deformation to the normal cone operations while keeping the local intersection number under concern invariant. 

\begin{theorem}[Proposition \ref{defo} and Proposition \ref{invar}]\label{m1}
Suppose $f$ has no blowup in codimension $0$. $[T^*_fM]$ is equivalent to a conical Lagrangian cycle $L$ in $T^*(M\times N)$ in the sense that both are members of an analytic family of cycles in $T^*(M\times N)$. In particular, $[T^*_fM]$ is homologous to $L$. The constructible function on $M\times N$ corresponding to $L$ is $(-1)^n\textup{Eu}^{\vee}_{f}$. 
\end{theorem}

The no blowup in codimension $0$ condition (see Definition \ref{noblow}) is weaker than the no blowup condition mentioned above. The weaker condition is enough to guarantee the existence of the deformation and the invariance of the local intersection number. Technically speaking we can always define the constructible function $\textup{Eu}^\vee_f$ under the weaker condition, but we need the stronger condition to attach a topological meaning to $\textup{Eu}^\vee_f$.

\begin{example}\label{e1}\cite[Theorem~4.5]{MR804052}
when $\dim N=1$, there exists a locally finite family (possibly empty family) of irreducible analytic subsets $Z_i \subset Z$ and positive integers $m_i$, such that $f$ is constant on each $Z_i$ and 
\begin{align*}
\textup{Eu}^{\vee}_f &= -\textup{Eu}^{\vee}_Z + \sum_i m_i \textup{Eu}^{\vee}_{Z_i} \\
-[T^*_fM] &\sim -[T^*_Z(M\times N)] + \sum_i m_i [T^*_{Z_i}(M\times N)], 
\end{align*}
where $\sim$ means the cycles on the two sides can be included in a family of analytic cycles in $T^*(M\times N)$ parametrised by $\mathbb{P}^1$. The analytic family is generated by deformation to the normal cone. The equivalence is the key ingredient to prove the holomorphic intersection formula in \cite{MR4749159}.
\end{example}

\vspace{1cm}
{\noindent\bf{Achknowledgement.}}
I want to express my deepest gratitude to Toru Ohmoto. In 2015 I was a postdoc fellow in KIAS. Uncertainty about the future drove me into a deep depression. One day near the end of 2015, I checked my email—a completely normal activity that the depression had prevented me from doing for several months—and found an email from Ohmoto inviting me to a conference. Though the conference was already over because of my delay, he gave me another invitation to visit him. While I was with him in Hokkaido in March 2016, he gave me his personal manuscript containing the main theorem as a conjectured formula and a topological proof under stronger conditions, encouraged me to find an algebraic proof. Some parts of the materials in that manuscript is presented now in \cite[\S 4.2.2]{MR5040859}. This formula motivated my research in the last decade. In 2022 I wrote a manuscript which contained most mathematical content presented in this paper, but I kept it private because I was stuck on the proofs of two technical lemmas. The gap is filled by ChatGPT 6.0, which provides the proof of Lemma \ref{lemma:invarianceintersection} and Lemma \ref{lemma:isolation}. I also want to give thanks to J\"{o}rg Schurmann and Claude Sabbah for their valuable comments and encouragement at various stages of this work. The author is supported by Chinese National Science Foundation of China (Grant No.11901214).

\section{perliminaries}

\subsection{The no blowup condition}

We collect necessary background materials about the no blowup condition. Basic reference for this subsection are \cite{MR760677,MR737936,MR804052}.

Let $M$ be a complex manifold, let $S$ be a reduced analytic space and let $Z$ be a closed analytic subset of $M\times S$. The projection $\pi: M\times S \to S$ induces a map $f: Z\to S$.
\begin{equation*}
\begin{tikzcd}
Z \arrow[rd,"f"] \arrow[r] & M \times S \arrow[d,"\pi"] \\
                                  & S
\end{tikzcd}
\end{equation*}
We always assume that there exists a Zariski open and dense subset $Z^\circ$ of $Z$ such that $f\vert_{Z^\circ}$ is smooth with $d$-dimensional fibres. By the upper semicontinuity of the fibre dimension, we know $\dim f^{-1}f(x)\geq d$ for any $x\in Z$. We define the relative conormal space $T^*_fM$ as the closure of the fibrewise conormal bundle $T^*_{f\vert_{Z^\circ}}M$ in $T^*M\times S$. The set $T^*_fM$ is closed and analytic in $T^*M\times S$ and it is irreducible if $Z$ is irreducible. If $S$ is a point, the set is simply denoted by $T^*_ZM$. Let $\tau_f: T^*_fM \to Z$ be the natural projection induced by $T^*M \times S \to M \times S$.

\begin{definition}[{\cite[Definition~4.1.1]{MR760677}}]\label{noblow}
$f$ has no blowup in codimension 0 if all nonempty fibres of the composition $f\circ\tau_f: T^*_fM \to S$ are of dimension $m$.
\end{definition}

Definition \ref{noblow} is a modification of the original one given in \cite{MR760677}. Indeed, the original definition requires that the map $f$ is surjective onto $S$. In our setting the image of $f$ is not necessarily an analytic subset of $S$. However, since the definition clearly implies that $\dim f^{-1}f(x)=d$ for all $x\in Z$, the local image of $f$ is a locally closed analytic subset of $S$ by the Remmert rank theorem \cite[p.295]{MR1131081}, so there's no harm to use this modified definition to include a slightly more general situation. Moreover, by upper semicontinuity of fibre dimension applied to $f\circ \tau_f$, the definition also implies that all fibres of $f\circ\tau_f$ have pure dimension $m$.

The no blowup in codimension 0 condition is stable under base change. Let $h: S' \to S$ be a holomorphic map. Consider the fibre product squares below where we suppose $f$ has no blowup in codimension 0.

\begin{equation*}
\begin{tikzcd}
T^*_fM \times_S S' \arrow[d,"\tau'_f"] \arrow[r] & T^*_fM \arrow[d,"\tau_f"] \\
Z \times _S S'\arrow[d,"f'"] \arrow[r] & Z \arrow[d,"f"] \\
               S'  \ar[r, "h"]                 & S
\end{tikzcd}
\end{equation*}

\begin{proposition}\cite[Corollary~4.2.1]{MR760677}\label{base-change}
Let $\Gamma$ be an irreducible component of $T^*_fM \times_S S'$. Then $\Gamma$ can be identified with the relative conormal space of the map $\tau'_f(\Gamma)\to S'$.
\end{proposition}

For the purpose of this paper, the most interesting case is that $S'$ is a singleton. Let $\Gamma$ be an irreducible component of $(f\circ\tau_f)^{-1}(s)$ where $s$ is an arbitrary point in $S$. The proposition implies that $\Gamma$ is the absolute conormal space of $\tau_f(\Gamma)$. Another way to show $(f\circ\tau_f)^{-1}(s)$ is conical Lagrangian is to apply the Lagrangian specialisation principle in \cite[Proposition 1.2.6]{MR966949}.

The no blowup in codimension 0 condition is also stable under deformation to the normal cone in the following sense. When $S$ is a complex manifold, $T^*_Z(M\times S)$ and $M \times T^*S$ are both analytic subsets of $T^*(M\times S)$. Let $\mathscr{C}(Z,S)$ be the normal cone to 
\begin{equation*}
\Big(M \times T^*S \Big)\cap T^*_Z(M\times S) 
\end{equation*}
in $T^*_Z(M\times S)$. 

\begin{proposition}\cite[Proposition~3.3]{MR804052}\label{induction}
Assume $f$ has no blowup in codimension 0 and $S$ is a complex manifold. Let $\mathscr{C}_i$ be an irreducible component of $\mathscr{C}(Z,S)$ and let $Z_i$ be its image in $Z$. Then $f\vert_{Z_i}$ has no blowup in codimension 0, and $\mathscr{C}_i \cap (T^*M \times S) = T^*_{f\vert_{Z_i}}M$. Moreover, at any $x\in Z_i$ if we shrink $Z_i$ so that its image by $f$ is a locally closed analytic subset $T_i$ of $S$,then $\mathscr{C}_i$ can be locally identified with the unique irreducible component of $T^*_{f\vert Z_i}M \times_{T_i} T^*_{T_i}S$ which dominates $Z_i$. 
\end{proposition}

\begin{remark}\label{remark:distinctness}
The last assertion of the proposition implies that the $Z_i$'s are a locally finite family of analytic subsets of $Z$ and distinct $\mathscr{C}_i$'s project to distinct $Z_i$'s.
\end{remark}

\begin{definition}[\cite{MR737936,MR476739}]\label{definition:noblowup}
We say the map $f:Z\to S$ has no blowup if $Z$ admits a stratification satisfying the Whitney conditions and the Thom $a_f$ condition between all pairs of adjacent strata.
\end{definition}

L\^{e} called such a stratification a good stratification in \cite[p.398]{MR476739}, but he considered only the case $S=\mathbb{C}$, and the terminology is not adopted by other authors. For example in \cite{MR760677} a good stratification means a stratification satisfying the conditions of Proposition \ref{prop:nbc0-stratification} below. We follow Sabbah's terminology in \cite{MR737936}. However, Sabbah demanded that
\begin{enumerate}
\item $f$ is stratified.
\item The stratification of $S$ is Whitney regular.
\item Pairs of adjacent strata of $Z$ satisfy the Thom $a_f$ condition.
\item Only those pairs of adjacent strata of $Z$ over a same stratum of $S$ satisfy Whitney conditions. 
\end{enumerate}
This is because Sabbah was considering a vanishing cycle theory locally on $S$ so he needed assumptions to apply the second isotopy lemma, while our concern is primarily on $Z$ so we need no knowledge of the stratification of $S$. In other words, we choose to graft L\^{e}'s condition to Sabbah's terminology. According to the following result, no blowup implies no blowup in codimension 0.

\begin{proposition}\cite[Theorem 7.1.1, equivalence (a) and (c)]{MR760677}\label{prop:nbc0-stratification}
The map $f$ has no blowup in codimension 0 if and only if $Z$ has a stratification such that $f$ has constant corank on each stratum and for any stratum $Z_i$ whose closure is an irreducible component of $Z$, any $Z_\alpha$ with $Z_\alpha \subset \overline{Z_i}$, the pair $(Z_i,Z_\alpha)$ satisfies the Thom $a_f$ condition.
\end{proposition}

\begin{remark}\label{remark:nb-stratification}
In either case, whether $f$ has no blowup or no blowup in codimension 0, the stratification on $Z$ induces a stratification on each fibre $f^{-1}(s)$ for any $s\in f(Z)$, say $f^{-1}(s)=\cup_\alpha S_\alpha$. Under the weaker no blowup in codimension 0 condition, we can already conclude that $|(f\circ\tau_f)^{-1}(s)|=\cup_\alpha T^*_{S_\alpha}M$. If $f$ has no blowup, we can conclude that $\{S_\alpha\}$ is a Whitney stratification. For any $z\in f^{-1}(s)$, we can find a small ball $B_\epsilon\subset M$ centred at $z$ such that $\partial B_\epsilon$ is transverse to each stratum $S_\alpha$. Then the Thom $a_f$ condition guarantees that $\partial B_\epsilon$ is also transverse to the stratifcation of $f^{-1}(s')$ for those $s'$ contained in a sufficiently small ball $B_\eta$ centred at $f(z)$.
\end{remark}

\begin{example}\label{example:maincases}
(i): A holomorphic map $f:Z\to \mathbb{C}$ such that the restriction of $f$ to any irreducible component of $Z$ is non-constant has no blowup accroding to \cite{MR476739}. (ii): The holomorphic map in the main theorem has no blowup. Indeed, any Whitney stratification of $M$ such that the critical space $C(f)$ is a union of strata and $M\setminus C(f)$ is a stratum satisfy the requirements of Definition \ref{definition:noblowup}. Because $f\vert_{C(f)}$ has discrete fibres, the fibres of $f\vert_{S_\alpha}$ are also discrete whenever $S_\alpha \subset C(f)$. Therefore the Thom $a_f$ condition is trivially satisfied.
\end{example}

\subsection{The local Euler obstruction and Chern class transformation}\label{section:CSM}
In the category of complex algebraic varieties (where the morphisms of the category are proper morphisms of algebraic varieties) we may consider two functors. They are the functor $F$ which associates the group of constructible functions to a variety, and the Borel-Moore homology functor $H_*$. The Chern class transformation was introduced by MacPherson \cite{MR361141} as a natural transformation $c_*: F\leadsto H_*$. Later on the theory was extended to algebraic varieties over an algebraically closed field of characteristic 0 in \cite{MR1063344}, and to complex analytic spaces which are embeddable in complex manifolds in \cite{MR1483328}. Below, we list a few definitions central to the calculus of constructible functions. General references for the theory of constructible functions and characteristic cycles are \cite{MR2143077,MS22}.

Let $M$ be an $m$-dimensional complex manifold.

\begin{definition}\label{def-lc}
A (conical) Lagrangian cycle is a summation
\begin{equation*}
\sum_i a_i[T^*_{Z_i}M]
\end{equation*}
with integer coefficients, where the $Z_i$'s are a locally finite family of closed irreducible analytic subsets of $M$. The summation is regarded as an element in $H_{2m}(\bigcup_iT^*_{Z_i}M)$.
\end{definition}

\begin{definition}
Let $Z$ be a complex analytic space. A constructible function $\alpha$ on $Z$ is an integer valued function $\alpha: Z \to \mathbb{Z}$ such that $\alpha^{-1}(k)$ is a constructible subset of $Z$ for any $k\in \mathbb{Z}$.
\end{definition}

Let $Z$ be an irreducible closed analytic subset of $M$. For any $z\in Z$, take a coordinate neighbourhood in $M$ centred at $z$ with coordinates $z_1, \ldots, z_m$. Let $r$ be the squared Euclidean distance function, i.e. $r(z) = |z_1|^2 + \ldots + |z_m|^2$. By Whitney (b) condition applied to a stratification of $Z$, 
\[
dr(B_\epsilon)\cap T^*_ZM =(z,0)\in T^*M
\]
whenever $0<\epsilon\ll 1$. Here $dr(B_\epsilon)$ is the graph of $dr$ over $B_\epsilon$ regarded as a subset of $T^*M$.

\begin{definition}\label{def-Eu}
The dual local Euler obstruction $\textup{Eu}^{\vee}_Z$ is the constructible function on $Z$ whose value at $z$ is the local intersection number
\begin{equation*}
\#_z\Big([dr(B_{\epsilon})] \cap [T^*_ZM]\Big).
\end{equation*} 
The value is independent of all the choices made above. The local Euler obstruction $\textup{Eu}_Z$ is defined by
\begin{equation*}
\textup{Eu}_Z = (-1)^{\dim Z}\textup{Eu}^{\vee}_Z. 
\end{equation*}
\end{definition}

\begin{remark}
$dr=2\sum_i (x_idx_i+y_idy_i)$ is a local section of the real cotangent bundle. The conormal space $T^*_ZM$ is regarded as a subset of the real cotangent bundle by forgetting the complex structure. The intersection product requires the specifications of orientations. The orientations for $B_{\epsilon}$, $T^*_ZM$ and $T^*M$ are all the canonical complex orientations. The detailed meaning of the topological intersection is explained in \cite[p.286]{MR2031639}.
\end{remark}

\begin{remark}\label{nonsingular}
$\textup{Eu}^\vee_Z(z) = (-1)^{\dim Z}$ if $z$ is a nonsingular point on $Z$.
\end{remark}

Any constructible function $\alpha$ on $M$ is uniquely expressible as a locally finite sum
\begin{equation*}
\alpha = \sum_i a_i \textup{Eu}^{\vee}_{Z_i}M
\end{equation*}
where $a_i \in \mathbb{Z}$ and the $Z_i$'s are a locally finite family of irreducible closed analytic subsets of $M$. Hence for each constructible function $\alpha$ we can associate the characteristic cycle
\begin{equation*}
\textup{CC}(\alpha) = \sum_i a_i [T^*_{Z_i}M].
\end{equation*} 
Conversely, for any Lagrangian cycle $L$ we can associate a constructible function $\textup{Eu}^{\vee}(L)$ whose value at $z$ is 
\begin{equation*}
\textup{Eu}^{\vee}(L) (z) = \#_z\Big([dr(B_{\epsilon})] \cap L\Big).
\end{equation*}
These two constructions are inverse to each other, i.e. $\textup{Eu}^\vee(\textup{CC}(\alpha))=\alpha$ for any constructible function $\alpha$.

Let $\pi_M: \mathbf{P}(T^*M) \to M$ be the projection and let  
\begin{equation*}
0 \to \xi_M \to  \pi^*_M T^*M \to \zeta_M \to 0
\end{equation*}
be the Euler sequence on $\mathbf{P}(T^*M)$ where $\xi_M$ is the tautological subbundle and $\zeta_M$ is the tautological quotient bundle. The Chern class transformation of $\alpha$, denoted by $c_*(\alpha)$, can be computed by the formula
\begin{equation*}
c_*(\alpha) = (-1)^{m} (\pi_M)_* \Big(c(\zeta^{\vee}_M) \cap \mathbf{P}(\textup{CC}(\alpha)\oplus 1)\Big).
\end{equation*}

\section{Characteristic cycles for relative dual Euler obstructions}\label{CC-rel-Eu}
\subsection{The relative dual local Euler obstruction} Let $f: Z \to N$ be a holomorphic map. Assume $f$ is open and has no blowup in codimension 0, $N$ is a complex manifold, and $f$ is factorised through $M\times N$. Take an arbitrary point $z\in Z$. As in the definition of the local Euler obstruction, we establish a local coordinate system centred at $z= (x, y) \in M\times N$ with coordinates $(z_1,\ldots,z_m, w_1,\ldots,w_n)$. Let $r_M = |z_1|^2 + \cdots +|z_m|^2$, $r_N = |w_1|^2+\cdots +|w_n|^2$, and $r=r_M + r_N$. Also let $B_{\epsilon}, B_{\eta}$ be sufficiently small open balls in $M$, $N$ centred at $x$, $y$ respectively. 

\begin{definition}\label{rel-Eu}
The relative dual local Euler obstruction is the function $\textup{Eu}^{\vee}_f: Z\to \mathbb{Z}$ given by
\begin{equation*}
\textup{Eu}^{\vee}_f(z) = (-1)^n\cdot \#_{(x,y)}\Big([dr(B_\epsilon \times B_\eta)] \cap [T^*_fM]\Big).
\end{equation*}
\end{definition}

Since both cycles appearing in the intersection product have real dimension $2(m+n)$, and the ambient space $T^*(M\times N)$ has real dimension $4(m+n)$, the local intersection number is well-defined if we can show for sufficiently small positive numbers $\epsilon,\eta$,
\[
dr(B_{\epsilon}\times B_{\eta}) \cap T^*_fM=\{z\}.
\] 
Let $(x',y',\alpha, \beta) \in T^*_fM$ be any point with $\alpha \in T^*_{x'}M$ and $\beta\in T^*_{y'}N$. By definition of $T^*_fM$, we necessarily have $\beta=0$. So if $(x',y',\alpha, \beta)$ also lies on $dr(B_\epsilon \times B_\eta)$, then necessarily $y'=0 \in B_\eta$. Therefore $(x,\alpha) \in (f\circ\tau_f)^{-1}(y)$, and the latter is a locally finite union of conormal spaces in $T^*M$ by Proposition \ref{base-change}. Hence by the paragraph before Definition \ref{def-Eu}, $x'=0$ and $\alpha=0$.

We may assume $f$ is surjective by replacing $N$ by $f(Z)$. By \cite[Corollary 4.3.13]{MR4294879} and the no blowup in codimension 0 assumption, the map $f\circ\tau_f$ defines an analytic family of cycles parametrised by $N$. Denote the cycle over $y\in N$ by $[T^*_fM]_y$. It is a conical Lagrangian cycle according to Proposition \ref{base-change}.

\begin{proposition}\label{prop:another-definition-reEu}
If $f: Z\to N$ is open and has no blowup in codimension 0, then
\begin{equation*}
\textup{Eu}^{\vee}_f(z) = \textup{Eu}^\vee([T^*_fM]_y)(x) 
\end{equation*}
\end{proposition}

\begin{proof}
By abuse of notation, we may regard $r_M$ as a function on $B_\epsilon$ and regard $r_N$ as a function on $B_{\eta}$. Therefore $dr_M:B_\epsilon \to T^*M$ is a local section of $T^*M$ and $dr_N:B_\eta \to T^*N$ is a local real section of $T^*N$. We need to show
\begin{equation*}
(-1)^n\cdot \#_{(x,y)}\Big([dr(B_\epsilon \times B_\eta)] \cap [T^*_fM]\Big) = \#_x \Big([dr_M(B_{\epsilon})] \cap [T^*_fM]_y \Big).
\end{equation*}
By the projection formula, we have
\begin{align*}
\#_{(x,y)}\Big([dr(B_\epsilon \times B_\eta)] \cap [T^*_fM]\Big) &= \#_{(x,y)}\Big([dr_M(B_\epsilon) \times dr_N(B_\eta)] \cap [T^*_fM]\Big) \\
&=\#_{(x,y)}\Big([dr_M(B_\epsilon) \times B_\eta] \cap (dr_N)^*[T^*_fM]\Big). 
\end{align*}

\begin{equation*}
\begin{tikzcd}
(f\circ \tau_f)^{-1}(y) \arrow[d] \arrow[r] & T^*_fM  \arrow[d]  \\
\{y\} \arrow[r] \arrow[d] & N \arrow[d] \\
B_\eta  \arrow[r,"dr_N"] & T^*N  
\end{tikzcd}
\end{equation*}
Consider the Cartesian diagram above. If we trivialise $T^*N$ locally by $d\omega_1,\cdots,d\omega_n$ and identify the holomorphic cotangent bundle with the real cotangent bundle by taking the real part of a complexified differential form, then $dr_N=\sum\bar{\omega}_id\omega_i$, hence 
\[
dr_N(B_\eta) \cap [N]=(-1)^n[y].
\]
Here the intersection number $(-1)^n$ is the determinant of the Jacobian matrix of the map $dr_N:\mathbb{R}^{2n} \to \mathbb{R}^{2n}$ written in complex variable by $(\omega_1,\cdots,\omega_n) \mapsto (\bar{\omega_1},\cdots,\bar{\omega_n})$. Therefore we have 
\[
(dr_N)^*[T^*_fM]=(-1)^n[T^*_fM]_y
\]
and the statement of the proposition follows.
\end{proof}

Now assume $f:Z\to N$ has no blowup. This implies that we can find $0<\eta \ll \epsilon \ll < 1$ such that $\partial B_\epsilon$ is transverse to the induced stratification of the fibre $f^{-1}(y')$ for any $y'\in B_{\eta}$ by Remark \ref{remark:nb-stratification}. If the stratification of $f^{-1}(y)$ is denoted by $\{T_\beta\}$, then
\[
dr_M(\partial B_\epsilon)\cap \Big(\bigcup_\beta T^*_{T_\beta}M\Big)=\emptyset.
\]

\begin{proposition}\label{prop:topologymeaning}
If $f: Z\to N$ is open and has no blowup, we have
\begin{align*}
\textup{Eu}^{\vee}_f(z) &=  \#_x \Big([dr(B_{\epsilon})] \cap [T^*_fM]_y \Big) \\
&= \int [dr(B_{\epsilon})] \cap [T^*_fM]_{y'} \\
&= \chi\Big(f^{-1}(y')\cap B_{\epsilon}, \textup{Eu}^\vee([T^*_fM]_{y'})\Big)
\end{align*}
\end{proposition}
\begin{proof}
The second equality follows from the invariance of intersection number under a small perturbation explained in Lemma \ref{lemma:invarianceintersection} below, and the third equality follows from the intersection formula \cite[Theorem 4.2]{MR2143077}. As a remark, if $Z^{\circ}$ is the open dense subset of $Z$ where $f\vert_{Z^\circ}$ is smooth and $y'\in B_\eta \cap f(Z^{\circ})$, then $[T^*_fM]_{y'}$ is the conormal space to $f^{-1}(y')$ in $M$ with multiplicity 1, so in this case $\textup{Eu}^\vee_f(z)$ is the Euler integral of the dual Euler obstruction of the nearby fibre $f^{-1}(y)\cap B_\epsilon$. 
\end{proof}

The invariance of the intersection number is based on the following general observation. It is well-known among experts, but we can't find a reference from the literature.

\begin{lemma}\label{lemma:invarianceintersection}
Let $W$ be a complex manifold of dimension $n$, let $X_s$ be a continuous family of analytic $k$-cycles in $W$ parametrised by a topological space $S$, and let $Y$ be a closed oriented $C^\infty$ submanifold of $W$ of real dimension $2n-2k$. Suppose there exists a compact subset $K\subset W$ such that $|X_s|\cap Y \subset K$ for every $s\in S$. Then the intersection number $Y \cdot X_s$ is locally constant on $S$.
\end{lemma}

\begin{proof}
Let $U$ be a relatively compact open subset of $W$ containing $K$. For any $s\in S$, let $X'_s=X_s\vert_U$, $Y'=Y\cap U$. By the commutative diagram
\[
\begin{tikzcd}[column sep=large, row sep=large]
H^{2k}(W,W\setminus Y)
  \arrow[r] \arrow[d]
  & H^{2k}(U,U\setminus Y') \arrow[d] \\
H^{2k}(|X_s|,|X_s|\setminus K)
  \arrow[r, "\cong"]
  & H^{2k}(|X'_s|,|X'_s|\setminus K)
\end{tikzcd},
\]
the Thom class of $Y$ in $W$ and the Thom class of $Y'$ in $U$ restrict to the same image via $H^{2k}(|X_s|,|X_s|\setminus K)\cong H^{2k}(|X'_s|,|X'_s|\setminus K)$. Hence $Y\cdot X_s=Y'\cdot X'_s$. Let $V$ be a relatively compact open neighbourhood of $K$ such that $\bar{V}\subset U$. There exists a tubular neighbourhood $T$ of $Y'$ in $U$ such that $\bar{T}\cap |X_s|\subset V$. Indeed, since 
\[
Y\cap |X_s|\cap (\bar{U}\setminus V) =\emptyset,
\]
and $|X_s|\cap(\bar{U}\setminus V)$ is compact, it has an open neighbourhood $A$ such that $\bar{A}\cap Y=\emptyset$. We can let $T$ be any tubular neighbourhood of $Y'$ disjoint from $\bar{A}$. Since $\bar{T}\cap(\bar{U}\setminus V)$ is also compact,
\[
|X_{s'}|\cap \bar{T}\cap (\bar{U}\setminus V) =\emptyset,
\]
for $s'$ in a neighbourhood of $s$ according to \cite[Proposition 4.2.9(1)]{MR4294879}. By \cite[Proposition 6.25]{MR658304}, we can find a differential $2k$-form $\varphi$ supported in $T$ representing the Thom class of $Y'$. Moreover, take a cut-off function $\chi\in C^\infty_c(U)$ such that $\chi\vert_V=1$. Then, for $s'$ in the same neighbourhood of $s$ we have
\[
Y'\cdot X'_{s'}=\int_{X'_{s'}}\varphi=\int_{X'_{s'}}\chi\varphi
\]
because $\textup{supp}(\varphi)\cap |X'_{s'}|\subset T\cap |X'_{s'}|\subset V$. For $s'$ in a possibly even smaller neighbourhood of $s$, the last integral is independent of $s'$ according to \cite[Theorem 4.2.44]{MR4294879}.
\end{proof}

\subsection{The deformation algorithm}
Under the assumption that $f:Z\to N$ is open and has no blowup in codimension 0, we will construct a deformation for $T^*_fM$ based on inductive applications of Proposition \ref{induction}. In the paragraph below, if two analytic cycles  $X$ and $Y$ of a complex manifold $W$ are members of one analytic family of cycles, we write $X\sim Y$.

\textit{Step} 0: There is an analytic family of cycles connecting $T^*_Z(M\times N)$ to $\mathscr{C}(Z,N)$ given by deformation to the normal cone. By Proposition \ref{induction}, the union of the irreducible components of $\mathscr{C}(Z,N)$ dominating irreducible components of $Z$ is $T^*_fM$. Therefore
\begin{equation}\label{s0}
[T^*_fM] \sim [T^*_Z(M\times N)] + \sum_{i\in I_0} m_{0i}[\mathscr{C}_{0i}]
\end{equation}
where the $\mathscr{C}_{0i}$'s are other irreducible components of $\mathscr{C}(Z,N)$ and the $m_{0i}$'s are negative integers. Let $Z_{0i}=\tau_f(\mathscr{C}_{0i})$. By Remark \ref{remark:distinctness} the $Z_{0i}$'s are distinct analytic subsets of $Z$ with $\dim Z_{0i} < \dim Z$, and they form a locally finite family. For any $z\in Z_{0i}$, replacing $Z_{0i}$ by a small neighbourhood of $z$ and letting $T_{0i} = f(Z_{0i})$ be the local image, $\mathscr{C}_{0i}$ can be identified with the irreducible component of $T^*_{f\vert_{Z_{0i}}}(M) \times_{T_{0i}} T^*_{T_{0i}}N$ dominating $Z_{0i}$.

\textit{Step} 1: For each $i\in I_0$, deformation to the normal cone gives an analytic family of cycles parametrised by $\mathbb{P}^1$ connecting $T^*_{Z_{0i}}(M\times N)$ to $\mathscr{C}(Z_{0i},N)$. Let $\mathscr{C}'_{0i}$ be the unique irreducible component of $\mathscr{C}(Z_{0i},N)$ donimating $Z_{0i}$ and let $\mathscr{C}_{0ij}$ be other irreducible components. Therefore
\[
[\mathscr{C}'_{0i}] \sim [T^*_{Z_{0i}}(M\times N)] + \sum_{j} n_{0ij} [\mathscr{C}_{0ij}]
\] 
for some negative integers $n_{0ij}$. By Proposition \ref{induction} $f\vert_{Z_{0i}}$ has no blowup in codimension 0, hence by the same proposition applied to $f\vert_{Z_{0i}}$, $\mathscr{C}'_{0i}$ is locally identified with the irreducible component of $T^*_{f\vert_{Z_{0i}}}(M) \times_{T_{0i}} T^*_{T_{0i}}N$ dominating $Z_{0i}$. Consequently, $\mathscr{C}'_{0i} = \mathscr{C}_{0i}$. Adding these families with multiplicities, we have
\begin{equation}\label{s1}
\sum_{i\in I_0}m_{0i}[\mathscr{C}_{0i}] \sim \sum_{i\in I_0}m_{0i}[T^*_{Z_{0i}}(M\times N)] + \sum_{i\in I_0}\sum_{j} m_{0i}n_{0ij} [\mathscr{C}_{0ij}].
\end{equation}
(Although \cite[p.376]{MR4294879} only explains addition of two families, one can easily generalise it to handle addition of infinitely many families, provided at each point of the parameter space the summands form a locally finite family of cycles. If each summand family is analytic, then the analyticity of the final summation family follows from \cite[Principle 4.3.5]{MR4294879}.)

Substituting \eqref{s1} into \eqref{s0} for each $i\in I_0$, combining terms and reindexing irreducible components $\mathscr{C}_{0ij}$'s, we have
\begin{equation}\label{s2}
[T^*_fM] \sim [T^*_Z(M\times N)] + \sum_{i\in I_0} m_{0i}[T^*_{Z_{0i}}(M\times N)] + \sum_{i\in I_{1}} m_{1i}[\mathscr{C}_{1i}] 
\end{equation}
where $m_{0i}<0$ for every $i\in I_0$ and $m_{1i}>0$ for every $i\in I_1$. We caution the reader that the analytic family underlying the equivalence \eqref{s2} is parametrised by the one point union of two copies of $\mathbb{P}^1$. Let $Z_{1i}$ be the image of $\mathscr{C}_{1i}$ in $Z$. The family $\{Z_{1i}\}_{i\in I_1}$ have the properties
\begin{enumerate}
\item $\{Z_{1i}\}_{i\in I_1}$ is locally finite.
\item $f\vert_{Z_{1i}}$ has no blowup in codimension 0 by Proposition \ref{induction}.
\item Each $Z_{1i}$ ($i\in I_1$) is a proper analytic subset of $Z_{0j}$ for some $j\in I_0$ by Proposition \ref{induction}, hence $\dim Z_{1i} < \dim Z -1$ for any $i\in I_1$.
\end{enumerate}
Theses properties allow us to repeatedly use the stategies above. For example, in the next step we deform $T^*_{Z_{1i}}(M\times N)$ to $\mathscr{C}(Z_{1i},N)$.

\textit{Step} $k$: we have 

\begin{equation*}
[T^*_fM] \sim [T^*_Z(M\times N)] + \sum^{k-1}_{j=0}\sum_{i\in I_j} m_{ji}[T^*_{Z_{ji}}(M\times N)] + \sum_{i\in I_{k}} m_{ki}[\mathscr{C}_{ki}] 
\end{equation*}
where $(-1)^{j+1}m_{ji}>0$ for any $0\leq j \leq k$ and $i\in I_j$. Let $Z_{ki}$ be the image of $\mathscr{C}_{ki}$ in $Z$. The family $\{Z_{ki}\}_{i\in I_k}$ satisfies analogous properties (i)-(iii) as above, except that $\dim Z_{ki} < \dim Z -k$.

The process will terminate at step $k+1$ if for any $i\in I_k$ the set $\mathscr{C}(Z_{ki},N)$ is irreducible. In other words, 
\[
\mathscr{C}(Z_{ki},N)=\mathscr{C}_{ki}=T^*_{f\vert_{Z_{ki}}}M; \quad \mathscr{C}_{ki}\sim T^*_{Z_{ki}}M.
\] 
This happens, for instance, if the map $f\vert_{Z_{ki}}$ is a stratified submersion onto a submanifold of $N$. It may take less than $\dim Z$ steps before the process terminates. For instance, when $\dim N=1$ the construction terminate at step 0. Summarising, and letting $I_{-1} = \{\bullet\}$ a singleton, $Z_{-1\bullet}=Z$, $m_{-1\bullet}=1$, we have

\begin{proposition}\label{defo}
There exists an integer $l$ with $-1 \leq l < \dim Z$, a non-empty index set $I_j$ and a locally finite family of irreducible analytic subsets $\{Z_{ji}\}_{i\in I_j}$ for any $-1\leq j \leq l$, such that 
\begin{equation}\label{decom}
[T^*_fM] \sim \sum^{l}_{j=-1}\sum_{i\in I_j} m_{ji}[T^*_{Z_{ji}}(M\times N)].  
\end{equation}
In addition, $f\vert_{Z_{ji}}$ has no blowup in codimension 0, $(-1)^{j+1}m_{ji}>0$, and $\dim Z_{ji}< \dim Z -j$ for any $i\in I_j$. The analytic family underlying the equivalence is parametrised by a chain of $l+1$ $\mathbb{P}^1$'s.
\end{proposition}

\subsection{The characteristic cycle of the relative obstruction}
Let $L$ be the Lagrangian cycle homologous to $[T^*_fM]$ appearing in Proposition \ref{defo}. For any $z=(x,y) \in Z$, let $r$, $B_\epsilon$, $B_\eta$ be those in Definition \ref{rel-Eu}. We will prove
\begin{proposition}\label{invar}
\begin{equation*}
\textup{Eu}^\vee_f(z) = (-1)^n \cdot \#_{z}\Big([dr(B_\epsilon \times B_\eta)] \cap L\Big).
\end{equation*}
\end{proposition}

\begin{corol}\label{coro:constructibility}
$\textup{Eu}^\vee_f$ is a constructible function.
\end{corol}

Since $L$ is constructed inductively, the following result implies Proposition \ref{invar}.
\begin{proposition}\label{prop:invariance-by-steps}
Let $f:Z\to N$ be a holomorphic map having no blowup in codimension 0 (we do not assume $f$ is open here). Then
\begin{equation*}
\#_{z}\Big([dr(B_\epsilon\times B_\eta)]\cap [T^*_Z(M\times N)]\Big) = \#_z \Big([dr(B_\epsilon \times B_\eta)] \cap [\mathscr{C}(Z,N)]\Big)
\end{equation*}
\end{proposition}

\begin{proof}
Consider the deformation from $T^*_Z(M\times N)$ to $\mathscr{C}(Z,N)$ and denote the total space of the deformation by $\Gamma \subset T^*(M\times N) \times \mathbb{P}^1$ (see \cite[Chapter 5]{MR1644323}). The projection to $\mathbb{P}^1$ gives a map $\Gamma \to \mathbb{P}^1$ such that $\Gamma_1 = T^*_Z(M\times N)$ and $\Gamma_{\infty} = \mathscr{C}(Z,N)$. Let $I=\mathbb{R}_{\geq 1}\cup \infty \subset \mathbb{P}^1_{\mathbb{C}}$. Restricting the $\mathbb{P}^1_{\mathbb{C}}$-family to $I$ gives a continuous family of cycles parametrised by $I$, which is still denoted by $\Gamma_t$. We will show in the next two lemmas that there exists sufficiently small positive numbers $\epsilon,\eta$ such that
\[
dr(B_\epsilon\times B_\eta)\cap |\Gamma_t|=\{z\}
\]
for any $t\in I$. The conclusion will follow from Lemma \ref{lemma:invarianceintersection}.
\end{proof}

\begin{lemma}\label{lemma:isolation}
There exists $\epsilon,\eta>0$ such that $|\Gamma_t| \cap dr(B_{\epsilon} \times B_{\eta}) = \{z\}$ for any $t\in \mathbb{R}_{>0}$.
\end{lemma}
\begin{proof}
Replace $M\times N$ by a product coordinate neighbourhood centred at $z=(0,0)$. Take an arbitrary $(x',y',\alpha,\beta) \in T^*(M\times N)$. By the construction of the deformation space $\Gamma$, $(x',y',\alpha,\beta) \in |\Gamma_1|$ if and only if $(x',y',t\alpha,\beta) \in |\Gamma_t|$. So $(x',y',dr_M(x'), dr_N(y')) \in |\Gamma_t|$ if and only if $(x',y',\frac{1}{t}dr_M(x')), dr_N(y')) \in |\Gamma_1|$. Using homogeneity, this is also equivalent to $(x',y',\frac{1}{1+t}dr_M(x')), \frac{t}{1+t}dr_N(y')) \in |\Gamma_1|$. Define a weighted distance function by $r_a=ar_M+(1-a)r_N$. Then 
\[
(dr)^{-1}(|\Gamma_t|)=(dr_{\frac{1}{1+t}})^{-1}(|\Gamma_1|).
\] 

Suppose no $\epsilon,\eta$ can be found. So we can find a sequence of point $p_i \in Z\setminus \{z\}$ converging to $z$ and a sequence of numbers $t_i>0$ such that $p_i\in (dr)^{-1}(|\Gamma_{t_i}|)$ for all $i$. Let $a_i=\frac{1}{1+t_i}\in (0,1)$, we have
\[
p_i\in (dr_{a_i})^{-1}(|\Gamma_1|), \quad \forall i.
\]

Take a subanalytic Whiteny stratification of $Z$. By local finiteness of the stratification and passing to a subsequence, we may assume there exists a stratum $S$ such that $p_i\in S$ for all $i$ and $\lim a_i=a_0\in [0,1]$. Since the set
\[
F=\{(p,a)\in (S\setminus z)\times (0,1) \ | \ p\in (dr)^{-1}_a(|\Gamma_1|)\}
\]
is subanalytic in $Z\times \mathbb{R}$ and $(z,a_0)\in \bar{F}$, by the curve selection lemma we can choose a real analytic curve $(x(s),y(s),a(s))\in F$ starting from $(z,a_0)$. Since $(x(s),y(s))\in S$ for all $s$, the tangent vector $(x'(s),y'(s))$ is tangent to $S$ hence is annihilated by any conormal vector to $Z$ according to Whitney (a) condition. In particular
\begin{equation}\label{contra}
0=\textup{Re}\Big(dr_{a(s)}(x(s),y(s))\cdot(x'(s),y'(s))\Big)=\frac{a(s)}{2}\frac{dr_M(t)}{ds}+\frac{1-a(s)}{2}\frac{dr_N(s)}{ds},
\end{equation}
where $r_M(s)=\sum^m_{j=1} \|x_j(s)\|^2$ and $r_N(s)=\sum^n_{j=1} \|y_j(s)\|^2$. This is a contradiction. Indeed, without loss of generality assume $x(s)\not\equiv0$ and write $x(s)=s^ku(s)$ with $k>0$ and $u(0)\neq 0$. Then 
\[
\frac{dr_M(s)}{ds}=\frac{d}{ds}\Big(\|u(s)\|^2s^{2k}\Big)=2k\|u(s)\|^2s^{2k-1}+O(s^{2k})>0
\]
for $0<s\ll 1$. Similarly $\frac{dr_N(s)}{ds}\geq 0$ (since the order of $y(s)$ may be 0, we do not have strict inequality here). On the other hand, the weights $\frac{a(s)}{2}$ and $\frac{1-a(s)}{2}$ are strictly positive when $s>0$, so the right hand side of the equation \eqref{contra} is strictly positive when $0<s\ll 1$.
\end{proof}

\begin{lemma}
There exists positive numbers $\epsilon_\infty$ and $\eta_\infty$ such that $|\Gamma_\infty|  \cap dr(B_{\epsilon_\infty} \times B_{\eta_\infty}) = \{z\}$.
\end{lemma}
\begin{proof}
replacing $M\times N$ by a product coordinate neighbourhood centred at $z=(0,0)$, we may assume $\mathscr{C}(Z,N)$ has only finite irreducible components $\mathscr{C}_1, \ldots, \mathscr{C}_s$. Let $Z_i$ and $T_i$ be the image of $\mathscr{C_i}$ in $M\times N$ and $N$ respectively. By Proposition \ref{induction} $\mathscr{C}_i \subset T^*_{f\vert Z_i}M \times_{T_i} T^*_{T_i}N$, so it is enough to prove $dr(B_{\epsilon_\infty}\times B_{\eta_\infty}) \cap T^*_{f\vert Z_i}M \times_{T_i} T^*_{T_i}N = (0,0,0,0)$ for $0< \epsilon_\infty \ll \epsilon$ and $0<\eta_\infty \ll \eta$. Suppose $(x',y',dr_M(x'),dr_N(y'))\in T^*_{f\vert Z_i}M \times_{T_i} T^*_{T_i}N$, then $(y',dr_N(y')) \in T^*_{T_i}N$. By the paragraph before Definition \ref{def-Eu}, we can choose a very small $\eta_\infty$ so that $y'=0$ whenever $(y',dr_N(y')) \in T^*_{T_i}N$ and $y' \in B_{\eta_\infty}$. Fix this choice $\eta_\infty$, so that we may assume $y'=0$ and $(x',dr_M(x'))\in |(T^*_{f\vert Z_i}M)_0|$ (the fibre over $0\in B_{\eta_\infty}$). By Proposition \ref{induction} $f\vert Z_i$ has no blowup in codimension 0, hence by Proposition \ref{base-change} $|(T^*_{f\vert Z_i}M)_0|$ is a union of conormal spaces. By the paragraph before Definition \ref{def-Eu} again we have $x'=0$, if we choose a very small $\epsilon_\infty$ and ask $x'\in B_{\epsilon_\infty}$.
\end{proof}

\begin{remark}
The restriction from the analytic family of cycles parametrised by $\mathbb{P}^1$ to the continuous family of cycles parametrised by $I$ is crucial for proving Proposition \ref{prop:invariance-by-steps}, because one cannot expect $dr(B_\epsilon\times B_\eta)\cap |\Gamma_t|$ is contained in a compact set for all $z\in \mathbb{P}^1$. For example, let $Z=\{y=x^2\}\subset \mathbb{C}^2$, $N=\mathbb{C}$ and $f:Z\to N$ is given by $(x,y)\mapsto y$. The map $dr:\mathbb{C}^2 \to T^*\mathbb{C}^2$ is given by $(x,y)\mapsto (x,y,\bar{x},\bar{y})$. A quick computation shows that $(x,y,\bar{x},\bar{y})\in \Gamma_t$ when $t=-\frac{1}{|x|^2}$. So the bad intersection points will enter into any small neighbourhood of $(0,0)$ when $t$ goes to infinity along $\mathbb{R}_{<0}\subset \mathbb{P}^1_{\mathbb{C}}$.
\end{remark}

\section{Chern classes for nearby and vanishing cycles}

In \S\ref{CC-rel-Eu}, we have thoroughly discussed $(-1)^n[T^*_fM]$. Rewrite the equivalence \eqref{decom} as
\begin{equation}\label{decomp}
(-1)^{\dim Z}[T^*_Z(M\times N)] \sim (-1)^{\dim Z}[T^*_fM] + (-1)^{\dim Z -1} \sum^{l}_{j=0}\sum_{i\in I_j} m_{ji}[T^*_{Z_{ji}}(M\times N)].
\end{equation}
Under the general assumption that $f$ has no blowup in codimension 0, this can be regarded as a virtual decomposition of $(-1)^{\dim Z}[T^*_Z(M\times N)]$ into the nearby and vanishing components. In this section we will further discuss this decomposition in the special cases already mentioned in Example \ref{example:maincases}:
\begin{enumerate}[(i)]
\item $\dim N =1$,
\item  $f: M\to N$ is a flat holomorphic map between two manifolds with relative dimension $l = m-n\geq 0$ whose fibres admit at most isolated singularities. We identify $M$ with $Z=gr(f) \subset M\times N$.
\end{enumerate} 
In both cases there are well-defined Milnor fibres, so $(-1)^n[T^*_fM]$ serves as a substitute for the characteristic cycle of the nearby cycle $\chi_f$ up to a sign by Theorem \ref{m1} and Proposition \ref{prop:topologymeaning}.

Case (i) is studied extensively in the literature, and the decomposition here is a rearrangement of the terms in Example \ref{e1}. In case (ii), $(-1)^{m}[T^*_Z(M\times N)]$ corresponds to the constant function $1_M$ by Remark \ref{nonsingular}, and $(-1)^{m}[T^*_fM]$ corresponds to the function $\chi_f=(-1)^l\textup{Eu}^\vee_f$ by Theorem \ref{m1} and Proposition \ref{prop:topologymeaning}. Since $1_M = \chi_f + (-1)^{l-1}\mu_f$ by definition of the Milnor number, comparing terms with equation \eqref{decomp} we see
\begin{equation}\label{CCmu}
(-1)^{m -1} \sum^{l}_{j=0}\sum_{i\in I_j} m_{ji}[T^*_{Z_{ji}}(M\times N)] = \textup{CC}\Big((-1)^{l-1}\mu_f\Big).
\end{equation}

Our next goal is to write down a general formula for $c_*(\textup{Eu}^\vee_f)$, and formulas for $c_*(\chi_f)$, $c_*(\mu_f)$ in case (ii). Let $\xi$ and $\zeta$ be the tautological subbundle and quotient bundle on $\mathbf{P}(T^*(M\times N)\vert_Z)$, and let $\pi: \mathbf{P}(T^*(M\times N)\vert_Z) \to Z$ be the projection.
\begin{corol}[to Theorem \ref{m1}]
\begin{equation*}
c_*(\textup{Eu}^\vee_f) = (-1)^{m+n-1}\pi_*\Big(c(\zeta^\vee)\cap \big((-1)^n\mathbf{P}(L)\big)\Big) = (-1)^{m-1}\pi_*\Big(c(\zeta^\vee) \cap [\mathbf{P}(T^*_fM)]\Big)
\end{equation*}
\end{corol} 
\begin{proof}
The first equality is standard. The $L$ in Theorem \ref{m1} is the Lagrangian cycle appearing in Proposition \ref{defo}. Since the deformation from $[T^*_fM]$ to $L$ clearly respect the $\mathbb{C}^*$ action on $T^*(M\times N)$, $[\mathbf{P}(T^*_fM)]$ and $\mathbf{P}(L)$ can also be put in an analytic family, hence they are also homologous \cite[p.420]{MR4294879}. The second equality follows because it is an equality of homology classes, not an equality of cycles.
\end{proof}

From now on, consider the special case (ii). To facilitate the discussion below, we introduce more notations and make some comments.
\begin{itemize}
\item $\xi_M$ and $\zeta_M$ are the tautological subbundle and quotient bundle on $\mathbf{P}(T^*M)$, $\xi_N$ and $\zeta_N$ are the tautological subbundle and quotient bundle on $\mathbf{P}(f^*T^*N)$, $\xi$ and $\zeta$ are the tautological subbundle and quotient bundle on $\mathbf{P}(T^*M \oplus f^*T^*N)$.
\item $\pi_M$, $\pi_N$ and $\pi$ are projections to $M$ from $\mathbf{P}(T^*M)$, $\mathbf{P}(f^*T^*N)$ and $\mathbf{P}(T^*M \oplus f^*T^*N)$.
\item $\mathbf{P}(T^*(M\times N)\vert_Z)$ can be identified with $\mathbf{P}(T^*M\oplus f^*T^*N)$. The conormal space $T^*_Z(M\times N)$ is identified with the graph of $-df: f^*T^*N \to T^*M$.
\item Let $C(f)$ be the critical space of $f$. The canonical morphism
\begin{equation*}
df: f^*T^*N \to T^*M
\end{equation*} 
restricts to a vector bundle injection on $M\setminus C(f)$. Let $\overline{(f^*T^*N)^{\circ}}$ be the closure of 
\begin{equation*}
df\Big(f^*T^*N \Big\vert_{M\setminus C(f)}\Big)
\end{equation*}
in $T^*M$ and denote its projectivisation in $\mathbf{P}(T^*M)$ by $\overline{\mathbf{P}(f^*T^*N)^{\circ}}$. It can be seen that $\overline{(f^*T^*N)^{\circ}}$ coincides with the relative conormal space $T^*_fM$.
\item The deformation from $T^*_Z(M\times N)$ to $\mathscr{C}(Z,N)$ is identified with deforming the family of graphs $\Gamma_t = gr(tdf) \subset T^*M \oplus f^*T^*N$ so that $\Gamma_1 = T^*_Z(M\times N)$ and $\Gamma_\infty = \mathscr{C}(Z,N)$.
\item On $\mathbf{P}(f^*T^*N)$, denote the composition $\xi_N \to \pi^*_Nf^*T^*N \xrightarrow{df} \pi^*_NT^*M$ by $\widetilde{df}$. We denote the family of graphs $gr(t\widetilde{df})\subset \pi^*_NT^*M \oplus \xi_N$ by $\gamma_t$, and this family arises as deformation from $\gamma_1$ to the normal cone $C_{\gamma_1 \cap \xi}\gamma_1$.
\item $\eta$ and $\tau$ are tautological subbundle and quotient bundle on $\mathbf{P}(\pi^*_NT^*M \oplus \xi_N)$. Note that $\eta$ is also the restriction of $\xi$ to $\mathbf{P}(\pi^*_NT^*M \oplus \xi_N)$.
\end{itemize}

\begin{corol}\label{c1}
\begin{equation*}
c_*(\chi_f)=(-1)^{n-1} c(f^*TN)\cap (\pi_M)_* \Big(c(\zeta^{\vee}_M)\cap \overline{\mathbf{P}(f^*T^*N)^{\circ}}\Big)
\end{equation*}
\end{corol}

\begin{proof}
Since $\mathbf{P}(T^*_fM) = \overline{\mathbf{P}(f^*T^*N)^{\circ}}$ is contained entirely in $\mathbf{P}(T^*M)$, we have
\begin{equation*}
\zeta \Big\vert_{\overline{\mathbf{P}(f^*T^*N)^{\circ}}} \cong \zeta_M \oplus \pi^*_MT^*N,
\end{equation*}
and
\begin{align*}
c_*(\chi_f) &= (-1)^{m-n}c_*(\textup{Eu}^\vee_f) \\
&= (-1)^{n-1}\pi_*\Big(c(\zeta^\vee_M \oplus \pi^*_Mf^*TN) \cap [\overline{\mathbf{P}(f^*T^*N)^{\circ}}]\Big) \\
&= (-1)^{n-1}(\pi_M)_*\Big(c(\zeta^\vee_M \oplus \pi^*_Mf^*TN) \cap [\overline{\mathbf{P}(f^*T^*N)^{\circ}}]\Big) \\
&=(-1)^{n-1} c(f^*TN)\cap (\pi_M)_* \Big(c(\zeta^{\vee}_M)\cap \overline{\mathbf{P}(f^*T^*N)^{\circ}}\Big).
\end{align*}
\end{proof}

\begin{remark}
The formula of Corollary \ref{c1} also works in case (i) with the additional assumption that $Z= gr(f)$ for non-constant holomorphic mappings $f: M\to N$ between manifolds. Indeed, our argument in the proof above only uses the facts (1) $Z$ is such a graph and (2) existence of local Milnor fibration at every point of $M$. In this context, $\overline{\mathbf{P}(f^*T^*N)^{\circ}}$ can be identified with $\textup{Bl}_{C(f)}M$ which is naturally embedded in $\mathbf{P}(T^*M)$, and $C(f)$ is sometimes called Jacobian subscheme in the literature. The formula can be expressed in terms of the dual Segre class of $C(f)$ in $M$. Let $h=c_1(\xi^\vee_M)$.
\begin{align*}
c_*(\chi_f) &= c(TM \oplus f^*TN) \cap (\pi_M)_*\Big( (1+h)^{-1} \cap [\textup{Bl}_{C(f)}M]\Big) \\
&= c(TM\oplus f^*TN) \cap (\pi_M)_* \Big([\textup{Bl}_{C(f)}M] - (1-h+h^2-\ldots)\cap [C_{C(f)}M] \Big) \\
&= c(TM\oplus f^*TN) \cap \Big([M] + (-1)^m \breve{s}(C(f),M) \Big).
\end{align*}  
This computation should be compared with the hypersurface Chern-Schwartz-MacPherson class formula obtained in \cite[Theorem 1.4]{MR1697199} and \cite[Corollary 2.4]{MR1795550}. Suppose $N= B_\eta \subset \mathbb{C}$. Note that the function $\chi$ and $\mu$ in \cite[\S 2]{MR1795550} are constructible functions supported on the hypersurface $f^{-1}(0)$, while our $\chi_f$ is supported on the ambient space $M$. So $\mu_f=(-1)^{m-1}(\chi_f -1)$ (in our notation) is equal to $\mu$ (in the notation of \cite{MR1795550}) but $\chi_f \neq \chi$. Taking into account of the notation difference, it is not difficult to see that when $N=B_\eta$ the computation above recovers \cite[Corollary 2.4]{MR1795550}.
\end{remark}

To calculate $c_*(\mu)$, the standard process which associates the Chern class to the characteristic cycle gives
\begin{equation*}
c_*(\mu_f) = (-1)^{(m+n-1)} \pi_*\Big(c(\zeta^\vee)\cap [\mathbf{P}(\textup{CC}(\mu_f))]\Big)
\end{equation*}
and $\textup{CC}(\mu_f)$ is given by the rather complicated expression \eqref{CCmu}. However, to get $c_*(\mu_f)$ we can replace the cycle $[\mathbf{P}(\textup{CC}(\mu_f))]$ by any cycle homologous to it in $H_{2(m+n-1)}(\mathbf{P}(T^*(M\times N)\vert_Z))$ in the formula above. Rearranging terms in \eqref{s0}, we get
\begin{equation*}
(-1)^{\dim Z}[T^*_Z(M\times N)] \sim (-1)^{\dim Z}[T^*_fM]+ (-1)^{\dim Z -1}\sum_{i\in I_0} m_{0i}[\mathscr{C}_{0i}].
\end{equation*}
Comparing with \eqref{decomp} and \eqref{CCmu}, we get
\begin{equation*}
\textup{CC}(\mu) \sim (-1)^{n-1} \sum_{i\in I_0} -m_{0i}[\mathscr{C}_{0i}].
\end{equation*}
Therefore
\begin{equation}\label{mue1}
c_*(\mu)= (-1)^{m} \pi_*\Big(c(\zeta^\vee)\cap \sum_{i\in I_0}-m_{0i}[\mathbf{P}(\mathscr{C}_{0i})]\Big).
\end{equation}
Note that $-\sum_{i\in I_0} m_{0i}[\mathscr{C}_{0i}]$ is the cycle complementary to $[T^*_fM]$ in $[\mathscr{C}(Z,N)]$ (recall $m_{0i}<0$), so that we only need deformation to the normal cone once.

We explain below an even better formula for $c_*(\mu)$ without reference to the graph embedding $gr(f) \subset M\times N$. For this, we need
\begin{proposition}
There is a proper modification $\kappa:\mathbf{P}(\pi^*_NT^*M \oplus \xi_N) \to \mathbf{P}(T^*M \oplus f^*T^*N)$ with the following properties.
\begin{enumerate}[(i)]
\item $\kappa$ respects deformation to the normal cone. In other words, if $\gamma \subset \mathbf{P}(\pi^*_NT^*M \oplus \xi_N) \times \mathbb{P}^1 $ is the total space of deformation for the family $\gamma_t$ and $\Gamma \subset \mathbf{P}(T^*M \oplus f^*T^*N) \times \mathbb{P}^1$ is so for $\Gamma_t$, then $(\kappa \times id_{\mathbb{P}^1}) \vert_\gamma : \gamma \to \Gamma$ is a proper modification. In particular, $\kappa_*([\mathbf{P}(\gamma_t)]) = [\mathbf{P}(\Gamma_t)]$ for any $t\in \mathbb{P}_1$. 
\item $\kappa^*\xi = \eta$.
\end{enumerate}
\end{proposition}

\begin{proof}
A point on $\mathbf{P}(\pi^*_NT^*M \oplus \xi_N)$ consists of the triple
\begin{itemize} 
\item a point $z\in M$, 
\item a 1-dimensional vector space $l_1 \subset f^*T^*N(z)$,
\item another 1-dimensional vector space $l_2 \subset T^*M(z) \oplus l_1$.
\end{itemize}
The proper modification $\kappa$ is given by $(z,l_1,l_2) \mapsto (z,l_2)$. Most statements are clear from this set-theoretic description of $\kappa$. The statement $\kappa_*([\mathbf{P}(\gamma_t)]) = [\mathbf{P}(\Gamma_t)]$ follows from the commutativity between proper push-forward and Gysin pull-back \cite[Theorem~6.2]{MR1644323}.
\end{proof}

The situation is best summarised in the following diagram.
\begin{equation}\label{lastdiagram}
\begin{tikzcd}
\mathbf{P}(\pi^*_NT^*M \oplus \xi_N) \arrow[r,"\kappa"] \arrow[d,"a"]& \mathbf{P}(T^*M \oplus f^*T^*N) \arrow[d,"\pi"] \\
\mathbf{P}(f^*T^*N) \arrow[r,"\pi_N"] & M
\end{tikzcd}
\end{equation}

We have 
\begin{equation*}
[\mathbf{P}(\gamma_\infty)] = [\mathbf{P}(C_{\gamma_1 \cap \xi}\gamma_1)] = [\mathbf{P}(\tilde{C} \oplus \xi_N\big\vert_{Z(\widetilde{df})})] + [\textup{Bl}_{Z(\widetilde{df})}\mathbf{P}(f^*T^*N)] 
\end{equation*}
by \cite[Example~18.1.6(d)]{MR1644323}. Here $Z(\widetilde{df}) \subset \mathbf{P}(f^*T^*N)$ is the analytic subspace defined by the zeros of $\widetilde{df}$, or equivalently the ideal sheaf $\mathscr{I}_{Z(\widetilde{df})}$ is the image of $\pi^*_NTM\otimes \xi_N \to \mathscr{O}_{\mathbf{P}(f^*T^*N)}$. The cone $\tilde{C}$ is defined by
\begin{equation*}
\tilde{C}:= \textup{Spec}\Big(\bigoplus_{i\geq 0} \big(\mathscr{I}^i_{Z(\widetilde{df})}\slash \mathscr{I}^{i+1}_{Z(\widetilde{df})}\big)\otimes \xi^{-i}_N\Big) 
\end{equation*}
 which is naturally a subcone of $\pi^*_NT^*M\big\vert_{Z(\widetilde{df})}$, and is {\em{not}} the normal cone to $Z(\widetilde{df})$ in $\mathbf{P}(f^*T^*N)$ which is defined by
 \begin{equation*}
 C:= C_{Z(\widetilde{df})}\mathbf{P}(f^*T^*N) = \textup{Spec}\Big(\bigoplus_{i\geq 0} \mathscr{I}^i_{Z(\widetilde{df})}\slash \mathscr{I}^{i+1}_{Z(\widetilde{df})} \Big).
 \end{equation*}
 
 However, be aware that 
 \begin{equation*}
 \mathbf{P}(\tilde{C} \oplus \xi_N\big\vert_{Z(\widetilde{df})}) = \mathbf{P}(C\oplus 1).
 \end{equation*}
 
 By the proposition above, we have
 \begin{align*}
 \kappa_*([\mathbf{P}(\gamma_\infty)]) &= [\mathbf{P}(\Gamma_\infty)] \\
 &= [\mathbf{P}(\mathscr{C}(Z,N))] \\
 &= [\overline{\mathbf{P}(f^*T^*N)^{\circ}}] + \sum_{i\in I_0}-m_{0i}[\mathbf{P}(\mathscr{C}_{0i})].
 \end{align*}
 
Since both $\textup{Bl}_{Z(\widetilde{df})}\mathbf{P}(f^*T^*N)$ and $\overline{\mathbf{P}(f^*T^*N)^{\circ}}$ are birational to $\mathbf{P}(f^*T^*N)$, we have $\kappa_*[\textup{Bl}_Z\mathbf{P}(f^*T^*N)] = [\overline{\mathbf{P}(f^*T^*N)^{\circ}}]$. Consequently,
\begin{equation}\label{pushC}
\kappa_* [\mathbf{P}(C \oplus 1)] = \sum_{i\in I_0}-m_{0i}[\mathbf{P}(\mathscr{C}_{0i})].
\end{equation}

\begin{proposition}\label{reso}
$\dim Z(\widetilde{df}) = n-1$, so that $\widetilde{df}$ is a regular section and  
\begin{equation*}
C=f^*T^*M \otimes \zeta^\vee_N \Big\vert_{Z(\widetilde{df})}.
\end{equation*}
\end{proposition}
\begin{proof}
We first construct a decomposition of the critical locus $C(f)$ into a locally finite disjoint union of connected constructible submanifolds
\begin{equation*}
C(f) = \bigcup_{i\in I} V_i
\end{equation*}
such that $f\vert_{V_i}$ has constant rank, see also \cite[p.26]{MR0436203}. Let $\{W_j\}$ ($j\in J_1$) be the family of connected components of $C(f)^\circ$ (the nonsingular part of $C(f)$), and let $V_j$ ($j\in J_1$) be the dense open subset of $W_j$ where $f\vert_{W_j}$ has maximal rank. By \cite[Theorem~1 on p.265]{MR1131081}, $W_j\setminus V_j$ is constructible in $M$. Therefore the set
\begin{equation*}
C(f)_1 = \bigcup_{j\in J_1} (W_j\setminus V_j) \cup C(f)^*
\end{equation*}
is a close analytic subset of $M$ with $\dim C(f)_1 < \dim C(f)$ (here $C(f)^*$ stands for the singular locus of $C(f)$). Repeat this procedure using $C(f)_1$ in replacement of $C(f)$ and define $V_j$ ($j\in J_2, J_3, \ldots$) inductively. The construction will terminate in finite steps. Finally, letting $I = \cup J_k$, we obtain the desired stratification of $C(f)$.

Since $f\vert_{C(f)}$ has discrete fibres, $f\vert_{V_i}$ also has discrete fibres, which implies that $\textup{rank}\, f\vert_{V_i} = \dim V_i$. The set
\begin{equation*}
S_i = \{(x,\beta)\in V_i \times_N T^*N \ \vert \ \beta \text{ annihilates } df(T_xV_i) \}
\end{equation*}
therefore has dimension $n$. On the other hand, the set
\begin{equation*}
S= \{(x,\beta)\in C(f) \times_N T^*N \ \vert \ \beta \text{ annihilates } df(T_xM) \}
\end{equation*}
is contained in $\cup_{i\in I} S_i$, hence $\dim S \leq n$. Since $Z(\widetilde{df})$ is the projectivisation of $S$, $\dim Z(\widetilde{df}) \leq n-1$. Because $Z(\widetilde{df})$ is the zeros of a section of a rank $m$ vector bundle $\pi^*_NT^*M \otimes \xi^\vee_N$ over a $(m+n-1)$-dimensional base $\mathbf{P}(f^*T^*N)$, the dimension of every irreducible component of $Z(\widetilde{df})$ is no less than $n-1$. Hence $Z(\widetilde{df})$ has pure dimension $n-1$ and this implies that $Z(\widetilde{df})$ is Cohen-Macaulay by \cite[Theorem~17.3, 17.4]{MR1011461}. The conclusion of the proposition follows from \cite[Proposition~7.1]{MR1644323}.
\end{proof}

\begin{remark}
If $f$ is Thom transversal \cite[4.10]{MR3112993}, then the mapping germ $f: (M,x)\to (N,f(x))$ is stable at every point $x\in M$ \cite[Proposition~4.12]{MR3112993}, and $Z(\widetilde{df})$ is a resolution of singularity for $C(f)$. 
\end{remark}

\begin{theorem}\label{m3}
\begin{equation*}
c_*(\mu_f) = (\pi_N)_*\Big( c(\zeta^{\vee}_N)\cap [Z(\widetilde{df})] \Big)
\end{equation*}
\end{theorem}

\begin{proof}
It follows from equations \eqref{mue1} and \eqref{pushC} and diagram \eqref{lastdiagram} that
\begin{align*}
c_*(\mu) &=(-1)^m\pi_{N*}a_*\Big(\frac{c(a^*\pi^*_NTM)c(a^*\pi^*_NTN)}{c(\xi^{\vee})} \cap [\mathbf{P}(C \oplus 1)] \Big) \\
&=(-1)^m\pi_{N*}a_*\Big(\frac{c(a^*\pi^*_NTM)c(a^*\pi^*_NTN)}{c(\eta^{\vee})} \cap [\mathbf{P}(C \oplus 1)] \Big) \\
&=(-1)^m\pi_{N*}a_*\Big(\frac{c(a^*\pi^*_NTM)c(a^*\xi^{\vee}_N)}{c(\eta^{\vee})}\cdot\frac{c(a^*\pi^*_NTN)}{c(a^*\xi^{\vee}_N)} \cap [\mathbf{P}(C \oplus 1)]\Big) \\
 &=(-1)^m\pi_{N*}a_*\Big(c(\tau^{\vee})c(a^*\zeta^{\vee}_N) \cap [\mathbf{P}(C \oplus 1)]\Big).
\end{align*}

By Proposition \ref{reso},  
\begin{equation*}
\mathbf{P}(C\oplus 1) = \mathbf{P}(\pi^*_NT^*M \oplus \xi_N)\Big\vert_{Z(\widetilde{df})}.
\end{equation*}

On the other hand, by \cite[Example~3.3.4]{MR1644323} we have
\begin{equation*}
a_*\Big(c(\tau^{\vee})\cap[\mathbf{P}(\pi^*_NT^*M \oplus \xi_N)\Big\vert_{Z(\widetilde{df})}]\Big) = (-1)^m[Z(\widetilde{df})].
\end{equation*}

Putting these equations together and using projection formula, the theorem is proved.
\end{proof}

Using Theorem \ref{m3} and Proposition \ref{reso}, the proof of the main theorem becomes an easy exercise in intersection theory. To simplify notation, let $h=c_1(\xi^\vee_N)$, $E=TM$, $F=TN$, and we tacitly omit writing $f^*$ and $\pi^*_N$ when pulling back vector bundles.

\begin{align*}
c_*(\mu)&= (\pi_N)_*\Big(c(\zeta^\vee_N) \cap [Z(\widetilde{df})]\Big) \\
&= c(F)\cap (\pi_N)_*\Big( (1+h)^{-1}c_m(E^\vee \otimes \xi^\vee_N) \cap [\mathbf{P}(f^*T^*N)] \Big) \\
&= c(F)\cap (\pi_N)_* \Big( \sum^{\infty}_{i=0} (-1)^ih^i \sum^{m}_{j=0} c_{m-j}(E^\vee)h^{j} \cap [\mathbf{P}(f^*T^*N)] \Big) \\
&= c(F)\cap (\pi_N)_* \Big(\sum^{\infty}_{i=0} \sum_{k+j=m+i} (-1)^ic_k(E^\vee)h^j \cap [\mathbf{P}(f^*T^*N)] \Big) \\
&= c(F) \cap (\pi_N)_* \Big(\sum^{\infty}_{i=0} \sum_{k+j=m+i}(-1)^{i+k}c_{k}(E)h^j \cap [\mathbf{P}(f^*T^*N)] \Big) \\
&= c(F) \cap \Big(\sum^{\infty}_{i=0} \sum_{k+j=m+i}(-1)^{i+k}c_k(E)s_{j-(n-1)}(F^\vee) \cap [M] \Big) \\
&= \sum^{\infty}_{i=0} \sum_{k+j=m+i}(-1)^{l+1}c(F)c_k(E)s_{j-(n-1)}(F) \cap [M] \\
&= \sum^{\infty}_{i=0} (-1)^{l+1}c(F)c_{m+i-(n-1)}(E-F) \cap [M] \\
&= \sum^{\infty}_{i=0} (-1)^{l+1}c(F)c_{l+1+i}(E-F) \cap [M].
\end{align*}

\bibliographystyle{plain}
\bibliography{liao}

\end{document}